%% file: article.tex
\documentclass[11pt, a4paper]{article}
\input{preamble.tex}

\title{The \quivertools package for SageMath and Julia}
\author{Pieter Belmans \and Hans Franzen \and Gianni Petrella}
\date{\today}

\begin{document}

\maketitle

\begin{abstract}
  We introduce \quivertools,
  a new software package,
  available in both a SageMath and Julia version,
  to study quivers and their moduli spaces of representations.
  Its key features are the computation of general subdimension vectors,
  leading to canonical decompositions,
  and checking the existence of (semi)stable representations,
  as well as the enumeration of Harder--Narasimhan types
  and related calculations for Teleman quantization.
  Computations related to intersection theory on quiver moduli are also implemented.
\end{abstract}


\section{Introduction}
The representation theory of quivers
and the geometry of moduli spaces of quiver representations
are rich subjects,
surveyed in \cite{MR2484736,MR3882963,MR3966814}.
They admit very explicit calculations:
their origin lies close to linear algebra,
making an algorithmic approach viable.

Many aspects can be understood
from the properties of the infinite root system
defined by the \emph{Kac form} \eqref{equation:kac-form},
the symmetrization of the Euler bilinear form of the quiver,
as evidenced by the program initiated by Kac \cite{MR0718127,MR0557581,MR0677715,MR0607162}.
For example, the dimension vectors of indecomposable representations
of a quiver correspond to the positive roots of the associated root system.
Many other properties of quiver representations
and their moduli spaces
are likewise encoded in properties of the associated root system
and Kac--Moody Lie algebra,
as surveyed in \cite{MR3966814}.

For other results
describing properties of quiver representations
and their moduli spaces,
the essential ingredient is the \emph{Euler form} \eqref{equation:euler-form} itself,
see,
e.g.,
its many appearances in \cite{MR2484736}.
This Euler form will be the main tool
for the type of results
treated in this article,
whose goal is to describe some of the features of \quivertools,
an implementation of many algorithms
pertaining to quivers and their moduli spaces of representations.
The software is available both
as a library for SageMath
and as a Julia package, see \cite{quivertools}.

\paragraph{Related software}
\quivertools does not treat \emph{individual} quiver representations,
for which explicit algorithms exist too.
This is functionality provided by \texttt{QPA} \cite{QPA}.
Likewise,
\quivertools does not deal with enumerative invariants of quivers (with potential).
This functionality is provided by \texttt{CoulombHiggs}~\cite{CoulombHiggs}
and \texttt{msinvar}~\cite{msinvar}.

\paragraph{Structure of the article}
In \cref{section:quivers,section:quiver-moduli} we recall some basic definitions
and results about quivers, their representations, and their moduli spaces.
In \cref{section:general,section:harder-narasimhan,section:teleman-quantization,section:intersection-theory}
we subsequently discuss how important representation-theoretic or algebro-geometric questions
can be rephrased in terms of the Euler form of the quiver,
and thus allow for an explicit algorithmic approach,
implemented in \quivertools.

\filbreak
We illustrate some of its features
using the running example of the~$4$-Kronecker quiver
\begin{center}
  \begin{tikzpicture}
    \node (1) at (0,0)              {$\bullet$};
    \node (2) at (2,0)              {$\bullet$};

    \draw[->, bend left  = 30] (1) edge (2);
    \draw[->, bend left  = 10] (1) edge (2);
    \draw[->, bend right = 10] (1) edge (2);
    \draw[->, bend right = 30] (1) edge (2);
  \end{tikzpicture}
\end{center}
and dimension vector~$(2, 3)$.

\iftoggle{sagecode}{\iftrue}{\iffalse}
Each SageMath example requires \quivertools to be imported.
The installation is done by running
\mintinline{shell}{sage --pip install git+https://github.com/QuiverTools/QuiverTools.git@v1.1},
and the import is done by executing
\mintinline{sage}{from quiver import *}.
\fi

\iftoggle{juliacode}{\iftrue}{\iffalse}
Each Julia example requires \quivertools to be imported.
The installation is done by running
\mintinline{julia}{using Pkg; Pkg.add(url="https://github.com/QuiverTools/QuiverTools.jl", rev="3d81a07")}
inside the Julia REPL,
and the import is done by executing
\mintinline{julia}{using QuiverTools}.
\fi

Not all existing functionality is discussed.
For more features
one is referred to the documentation on \url{https://quiver.tools}.

\medskip

\paragraph{Acknowledgements}
P.B.~was partially supported by the Luxembourg National Research Fund (FNR--17113194).
H.F.~was partially supported by the Deutsche Forschungsgemeinschaft (DFG, German Research Foundation)
SFB-TRR~358/1~2023 ``Integral Structures in Geometry and Representation Theory'' (491392403).
G.P.~was supported by the Luxembourg National Research Fund (FNR--17953441).

\section{Representation theory of quivers}
\label{section:quivers}
A \emph{quiver} $Q$ is a directed multigraph,
given by a finite set of vertices~$Q_0$
and a finite set of arrows~$Q_1$.
Each arrow $\alpha \in Q_1$ has a \emph{source} $\source(\alpha)$
and a \emph{target} $\target(\alpha)$,
which are vertices in $Q_0$.

The \emph{Euler form} of a quiver $Q$ is the bilinear form
\begin{equation}
  \label{equation:euler-form}
  \langle\hyphen,\hyphen\rangle\colon
  \mathbb{Z}^{Q_0} \times \mathbb{Z}^{Q_0} \to \mathbb{Z}:
  (\tuple{d},\tuple{e}) \mapsto \sum_{i \in Q_0} d_i e_i - \sum_{a \in Q_1} d_{\source(a)} e_{\target(a)}.
\end{equation}
The motivation behind this definition and terminology
will be clear following \eqref{equation:euler-form-interpretation}.
By symmetrizing the Euler form
we obtain the so-called \emph{Kac form} on $\mathbb{Z}^{Q_0}$,
given by
\begin{equation}
  \label{equation:kac-form}
  (\mathbf{d},\mathbf{e})
  \colonequals
  \langle\mathbf{d},\mathbf{e}\rangle
  +
  \langle\mathbf{e},\mathbf{d}\rangle.
\end{equation}

We fix an algebraically closed field $\field$ of characteristic $0$.
A \emph{representation} of $Q$ is a collection of~$\field$\dash vector spaces~$V_i$
for every $i \in Q_0$,
and linear maps $V_{\alpha}\colon V_{\source(\alpha)} \to V_{\target(\alpha)}$
for every $\alpha \in Q_1$.
We say that a representation has \emph{dimension vector} $\tuple{d} \in \mathbb{Z}^{Q_0}$
if for all $i \in Q_0$, $\dim(V_i) = d_i$.
In this case we write~$\dimvec(V) = \tuple{d}$.

Given two representations $V$ and $W$ of $Q$,
a \emph{morphism} of representations from $V$ to $W$
is a collection of linear maps
$\{\phi_i\colon V_i \to W_i\}_{i\in Q_0}$,
such that for every arrow $\alpha \in Q_1$,
one has
\begin{equation}
  W_{\alpha} \circ \phi_{\source(\alpha)} = \phi_{\target(\alpha)} \circ V_{\alpha}.
\end{equation}

One of the first results on representations of quivers
is provided by Kac's theorem \cite{MR0557581}.
A non-zero representation~$V$ is said to be \emph{indecomposable}
if it cannot be written as the direct sum of two non-zero representations.
On the other hand,
Kac considers the symmetric bilinear form \eqref{equation:kac-form} of a quiver,
equipping~$\mathbb{Z}^{Q_0}$ with the structure of a root system,
for which the following can be shown.
\begin{theorem}[Kac]
  \label{theorem:kac}
  Let~$Q$ be a quiver.
  The following are equivalent for a dimension vector~$\mathbf{d}$:
  \begin{enumerate}
    \item there exists an indecomposable representation~$V$ for which~$\dimvec(V)=\mathbf{d}$;
    \item $\mathbf{d}$ is a root in the root system,
      i.e.,~$(\mathbf{d},\mathbf{d})\leq 2$.
  \end{enumerate}
\end{theorem}

The next two ingredients will form the basis for geometric considerations in what follows.
If we fix a dimension vector $\tuple{d}$, a representation corresponds
to a point in the \emph{representation space}
\begin{equation}
  \Rep(Q, \tuple{d}) \colonequals \prod_{\alpha \in Q_1}\Mat_{\target(\alpha), \source(\alpha)}(\field).
\end{equation}
The group $\GLd \colonequals \prod_{i \in Q_0} \GL_{d_i}(\field)$
acts on $\Rep(Q,\tuple{d})$ by base change;
its orbits are, by definition,
isomorphism classes of representations of $Q$.

\section{General representations of quivers}
\label{section:general}
The building block for many algorithms for moduli of quiver representations
is the notion of general representations,
introduced by Schofield \cite{MR1162487}.

\begin{definition}
  Let~$\mathcal{P}$ be a property of a representation of a quiver~$Q$.
  We say that the \emph{general representation} of dimension vector~$\tuple{d}$
  satisfies~$\mathcal{P}$
  if there exists a non-empty Zariski-open subset~$U$ of~$\Rep(Q,\tuple{d})$
  such that every representation in~$U$ satisfies~$\mathcal{P}$.
\end{definition}

The representations of a quiver $Q$
form an abelian category whose objects have all finite length,
and thus a Krull--Schmidt category.
Accordingly, each representation can be decomposed uniquely into a finite direct sum
of indecomposable ones.
The shape of this decomposition
for the general representation in $\Rep(Q, \tuple{d})$
can be described in terms of
the quiver $Q$ and the dimension vector $\tuple{d}$.
This is closely related to the following class of dimension vectors.

\begin{definition}
  \label{definition:general-subdimension-vector}
  A dimension vector $\tuple{e} \leq \tuple{d}$
  is called a \emph{general subdimension vector} of $\tuple{d}$
  if the general representation in $\Rep(Q,\tuple{d})$
  admits a subrepresentation of dimension vector $\tuple{e}$.
  We introduce the notation~$\tuple{e} \hookrightarrow \tuple{d}$
  to indicate that $\tuple{e}$ is a general subdimension vector of $\tuple{d}$.
\end{definition}

The subset of $\Rep(Q,\tuple{d})$ of representations admitting
a subrepresentation of fixed subdimension vector can be shown to be Zariski closed.
General subdimension vectors $\tuple{e}$ are thus characterised
by the (a priori) stronger requirement
that \emph{every} representation of dimension vector $\tuple{d}$ admits a
subrepresentation of dimension vector $\tuple{e}$ (see \cite[Theorem~3.3]{MR1162487}).
There exists a third and effective characterisation
of general subdimension vectors \cite[Theorem~5.2]{MR1162487}.
\begin{theorem}[Schofield]
  \label{theorem:general-subdimension-vector-characterization}
  Let $Q$ be a quiver and $\tuple{d}$ a dimension vector.
  A subdimension vector~$\tuple{e}\leq\tuple{d}$ is general
  if and only if,
  for every general subdimension vector $\tuple{e}'$ of $\tuple{e}$,
  the inequality $\langle \tuple{e}',\tuple{d} -\tuple{e}\rangle \geq 0$ holds.
\end{theorem}
The Euler form appearing in this recursive characterisation
is the first of many instances where it plays an essential role.
We describe an equivalent characterization,
allowing us to introduce concepts
which will be useful later.
\begin{definition}
  Let $Q$ be a quiver,
  and let~$\tuple{d}$ and $\tuple{d'}$ two dimension vectors.
  We define the \emph{general hom}~$\hom(\tuple{d},\tuple{d'})$
  and the \emph{general ext} $\ext(\tuple{d},\tuple{d'})$ by
  \begin{align}
    \hom(\tuple{d}, \tuple{d'})
    &\colonequals
    \min\left\{
      \dim(\Hom(V,W))\mid
      V \in \Rep(Q,\tuple{d}),
      W \in \Rep(Q,\tuple{d'})
    \right\}, \\
    \ext(\tuple{d}, \tuple{d'})
    &\colonequals
    \min\left\{
      \dim(\Ext(V,W))\mid
      V \in \Rep(Q,\tuple{d}),
      W \in \Rep(Q,\tuple{d'})
    \right\}.
  \end{align}
\end{definition}
Since the functions
\begin{equation}
  \begin{gathered}
    \dim(\Hom(-,-))\colon \Rep(Q,\tuple{d}) \times \Rep(Q,\tuple{d'}) \to \mathbb{Z} \\
    \dim(\Ext(-,-))\colon \Rep(Q,\tuple{d}) \times \Rep(Q,\tuple{d'}) \to \mathbb{Z}
  \end{gathered}
\end{equation}
are upper semi-continuous, they reach their minimum on an open subset of their domain.
The general representations $V \in \Rep(Q,\tuple{d})$ and $W \in \Rep(Q,\tuple{d'})$
thus satisfy $\dim(\Hom(V, W)) = \hom(\tuple{d},\tuple{d'})$
and $\dim(\Ext(V, W)) = \ext(\tuple{d},\tuple{d'})$.

The general $\hom$ and $\ext$ turn out to be closely related
to the Euler form.
Given two dimension vectors $\tuple{d}$ and $\tuple{d'}$,
and for every $V \in \Rep(Q,\tuple{d})$ and every $W \in \Rep(Q,\tuple{d'})$,
there exists a 4-term exact sequence
\begin{equation}
  \label{equation:4-term-sequence}
  0 \to
  \Hom(V, W) \to
  \bigoplus_{i\in Q_0}\Hom_\field(V_i,W_i) \to
  \bigoplus_{\alpha\in Q_1}\Hom_\field(V_{\source(\alpha)},W_{\target(\alpha)}) \to
  \Ext(V, W) \to
  0,
\end{equation}
where the middle morphism sends~$\{\phi_i\colon V_i\to W_i\}_{i\in Q_0}$
to~$\{W_\alpha\circ\phi_{\source(\alpha)}-\phi_{\target(\alpha)}\circ V_\alpha\}_{\alpha\in Q_1}$.
This 4-term sequence \eqref{equation:4-term-sequence}
(and its variations, see, e.g., \eqref{equation:4-term-global-sequence})
are an important tool in the study of the geometry of quiver moduli.

Computing the dimensions in the sequence \eqref{equation:4-term-sequence} tells us that
for every $V \in \Rep(Q,\tuple{d}), W \in \Rep(Q,\tuple{d'})$,
we have
\begin{equation}
  \label{equation:euler-form-interpretation}
  \langle \tuple{d}, \tuple{d'}\rangle =
  \dim(\Hom(V, W)) - \dim(\Ext(V,W)) =
  \hom(\tuple{d}, \tuple{d'}) - \ext(\tuple{d},\tuple{d'}).
\end{equation}
The following central result of Schofield
allows to compute general $\hom$ and $\ext$
using the Euler form and general subdimension vectors, see~\cite[Theorem~5.4]{MR1162487}.
We recall the following notation from op.~cit.,
dual to \cref{definition:general-subdimension-vector}:
for a dimension vector~$\mathbf{f}$
such that the general representation in~$\Rep(Q,\mathbf{d})$
admits a quotient representation of dimension vector~$\mathbf{f}$,
we write~$\mathbf{d}\twoheadrightarrow\mathbf{f}$.

\begin{theorem}[Schofield]
  \label{theorem:general-ext-computation}
  Let $\tuple{d}$ and $\tuple{d'}$ be dimension vectors for the quiver $Q$.
  Then,
  \begin{equation}
    \ext(\tuple{d},\tuple{d'})
    = \max_{\tuple{e} \hookrightarrow \tuple{d}}\{-\langle \tuple{e},\tuple{d'}\rangle\}
    = \max_{\tuple{d}' \twoheadrightarrow \tuple{f}}\{-\langle \tuple{d},\tuple{f}\rangle\}
    = \max_{\substack{\tuple{e} \hookrightarrow \tuple{d} \\ \tuple{d}' \twoheadrightarrow \tuple{f}}}\{-\langle \tuple{e},\tuple{f}\rangle\}.
  \end{equation}
\end{theorem}

In light of \cref{theorem:general-ext-computation},
\cref{theorem:general-subdimension-vector-characterization} can be
rephrased in terms of the general extension group.
\begin{corollary}
  A subdimension vector $\tuple{e}$ of $\tuple{d}$
  is general if and only if $\ext(\tuple{e}, \tuple{d - e}) = 0$.
\end{corollary}
We illustrate these notions with \quivertools.

\iftoggle{sagecode}{\iftrue}{\iffalse}
\begin{sagesnippet}
|\sage| Q = KroneckerQuiver(4); d = (2, 3);
|\sage| Q.all_general_subdimension_vectors(d)
[(0, 0), (0, 1), (0, 2), (0, 3), (1, 3), (2, 3)]
|\sage| Q.general_ext((1, 2), d)
4
|\sage| Q.general_hom((1, 2), d)
0
\end{sagesnippet}
\fi
\iftoggle{juliacode}{\iftrue}{\iffalse}
\begin{minted}{jlcon}
julia> Q, d = kronecker_quiver(4), [2, 3];

julia> all_general_subdimension_vectors(Q, d)
6-element Vector{StaticArraysCore.SVector{2, Int64}}:
[0, 0]
[0, 1]
[0, 2]
[0, 3]
[1, 3]
[2, 3]

julia> general_ext(Q, [1, 2], d)
4

julia> general_hom(Q, [1, 2], d)
0
\end{minted}
\fi

As stated above, the
dimension vectors of the indecomposable summands in the
decomposition of the general representation
of $Q$ of dimension vector $\tuple{d}$ only depend on $Q$ and $\tuple{d}$.
This is a result of Kac \cite{MR0557581}, which we state below.
\begin{theorem}[Kac]
  Let $Q$ be a quiver and $\tuple{d}$ be a dimension vector.
  There exists a unique collection of dimension vectors $\{\tuple{d}^i\}_{i = 1}^{n}$
  summing to $\tuple{d}$ and
  such that for all $i \neq j$,
  \begin{equation}
    \ext(\tuple{d}^i, \tuple{d}^j) = 0.
  \end{equation}
  Moreover, the general representation of $Q$ of dimension vector $\tuple{d}$
  decomposes as a direct sum $\bigoplus_{i=1}^nV^{i}$, where $\dimvec(V^{i}) = \tuple{d}^i$,
  and~$V^i$ is indecomposable.
\end{theorem}
This decomposition is accordingly called the \emph{canonical decomposition}
of the dimension vector~$\tuple{d}$.
In fact,
the dimension vectors $\{\tuple{d}^i\}_{i = 1}^{n}$
appearing in the canonical decomposition
are not just roots,
but they are \emph{Schur roots}:
the endomorphism algebra of a general representation
of any of these dimension vectors
is the base field.

An effective characterization of the canonical decomposition
follows from~\cref{theorem:general-ext-computation},
which can be used to recursively decompose a dimension vector
into general summands, as well as to determine which dimension vectors
are indecomposable.

We illustrate these notions with \quivertools,
with the 3-vertex example illustrating
what happens in \cite[Example~11.1.4]{MR3727119}.

\iftoggle{sagecode}{\iftrue}{\iffalse}
\begin{sagesnippet}
|\sage| Q1 = KroneckerQuiver(4);
|\sage| Q1.canonical_decomposition((2, 3))
((2, 3),)
|\sage| Q2 = ThreeVertexQuiver(1, 1, 1)
|\sage| Q2.is_root((1, 2, 1))
True
|\sage| Q2.is_schur_root((1, 2, 1))
False
|\sage| Q2.canonical_decomposition((1, 2, 1))
((0, 1, 0), (1, 1, 1))
\end{sagesnippet}
\fi
\iftoggle{juliacode}{\iftrue}{\iffalse}
\begin{minted}{jlcon}
julia> Q1 = kronecker_quiver(4);

julia> canonical_decomposition(Q1, [2, 3])
1-element Vector{Vector{Int64}}:
 [2, 3]

julia> Q2 = three_vertex_quiver(1, 1, 1);

julia> is_root(Q2, [1, 2, 1])
true

julia> is_schur_root(Q2, [1, 2, 1])
false

julia> canonical_decomposition(Q2, [1, 2, 1])
2-element Vector{StaticArraysCore.SVector{3, Int64}}:
 [0, 1, 0]
 [1, 1, 1]
\end{minted}
\fi

\section{Moduli spaces of representations of quivers}
\label{section:quiver-moduli}
As introduced in \cref{section:quivers},
the general linear group $\GLd$ acts on the representation space $\Rep(Q,\tuple{d})$
via change of basis.
One can consider the spectrum of the invariant ring:
\begin{equation}
  \modulispace(Q,\tuple{d}) \colonequals \Spec(\mathcal{O}(\Rep(Q,\tuple{d}))^{\GLd}).
\end{equation}
It is a moduli space, parametrising semisimple representations.
This invariant ring
is generated by the trace of products of matrices along
oriented cycles in the quiver \cite{MR958897}.
In particular, if $Q$ is acyclic then $\modulispace(Q,\tuple{d})$ is a point.

In order to obtain a better behaved moduli space,
the following notion of stability,
involving a choice of \emph{stability parameter} $\theta$,
was introduced in \cite{MR1315461}.
\begin{definition}
  Let $\theta \in \Hom(\mathbb{Z}^{Q_0}, \mathbb{Z})$ be
  such that $\theta(\tuple{d}) = 0$.
  A representation $V \in \Rep(Q,\tuple{d})$ is said to be \emph{$\theta$-stable},
  respectively \emph{$\theta$-semistable},
  if for all of its nonzero proper subrepresentations $W$, we have
  $\theta(\dimvec(W)) < 0$,
  respectively $\theta(\dimvec(W)) \leq 0$.
\end{definition}
We denote by $\Rep^{\theta\stable}(Q,\tuple{d})$
and $\Rep^{\theta\semistable}(Q,\tuple{d})$
the Zariski open subsets of $\Rep(Q,\tuple{d})$
of $\theta$-stable, respectively $\theta$-semistable
representations.
As it turns out, both
admit geometric quotients by the action of $\GLd$,
the former quotient is always smooth and
the latter is always
projective over $\Rep(Q,\tuple{d})/\GLd$.
These geometric quotients are the better behaved moduli spaces sought.
We denote them by
\begin{align}
  \modulispace^{\theta\stable}(Q,\tuple{d})
  &\colonequals \Rep^{\theta\stable}(Q,\tuple{d}) \gitquot_{\theta} \GLd, \\
  \modulispace^{\theta\semistable}(Q,\tuple{d})
  &\colonequals \Rep^{\theta\semistable}(Q,\tuple{d}) \gitquot_{\theta} \GLd.
\end{align}
For more on their construction and properties,
see, e.g., \cite{MR2484736}.
One standard choice of stability parameter is
the \emph{canonical} stability parameter for~$\mathbf{d}$,
defined as~$\theta_{\mathrm{can}}\colonequals\langle\mathbf{d},-\rangle-\langle-,\mathbf{d}\rangle$.
For~$Q$ acyclic this will often lead to smooth projective Fano varieties \cite{MR4352662},
making them particularly relevant for algebraic geometers.

Having introduced these moduli spaces,
we naturally wish to know for which pairs $(Q, \tuple{d})$
they are nonempty, i.e.,
for which such pairs (semi)stable representations exist.
These questions can be answered by verifying,
for a given stability parameter $\theta$,
whether the necessary (semi)stability inequalities hold for general subdimension vectors
(see, e.g.,~\cite[Theorem~3]{MR1758750}),
which makes for an explicit algorithm,
building upon \cref{section:general}.
We illustrate these notions with \quivertools:
in our running example,
the canonical stability condition
(and its rescallings, which do not affect stability)
is the only one for which the moduli space is interesting.

\iftoggle{sagecode}{\iftrue}{\iffalse}
\begin{sagesnippet}
|\sage| Q = KroneckerQuiver(4); d = (2, 3);
|\sage| theta = Q.canonical_stability_parameter(d)
|\sage| Q.has_stable_representation(d, theta)
True
|\sage| Q.has_stable_representation(d, -theta)
False
\end{sagesnippet}
\fi
\iftoggle{juliacode}{\iftrue}{\iffalse}
\begin{minted}{jlcon}
julia> Q, d = kronecker_quiver(4), [2, 3];

julia> theta = canonical_stability(Q, d);

julia> has_stables(Q, d, theta)
true

julia> has_stables(Q, d, -theta)
false
\end{minted}
\fi

As explained in \cite{MR2484736,MR2484736}
we can compute various geometric invariants
of quiver moduli spaces.
We compute these using \quivertools for our running example.

\iftoggle{sagecode}{\iftrue}{\iffalse}
\begin{sagesnippet}
|\sage| Q = KroneckerQuiver(4); d = (2, 3);
|\sage| theta = Q.canonical_stability(d);
|\sage| M = QuiverModuliSpace(Q, d, theta);
|\sage| M.dimension()
12
|\sage| M.picard_rank()
1
|\sage| M.index()
4
\end{sagesnippet}
\fi
\iftoggle{juliacode}{\iftrue}{\iffalse}
\begin{minted}{jlcon}
julia> Q, d, theta = kronecker_quiver(4), [2, 3];

julia> theta = canonical_stability(Q, d);

julia> M = QuiverModuliSpace(Q, d, theta);

julia> dimension(M)
12

julia> picard_rank(M)
1

julia> index(M)
4
\end{minted}
\fi

\section{Harder--Narasimhan stratification}
\label{section:harder-narasimhan}
For an arbitrary action of a reductive group $\mathrm{G}$ on an affine variety,
Hesselink has shown in \cite{MR514673}
that the complement of the semistable locus
admits a certain stratification into locally closed subsets.
This stratification is not unique,
depending on a choice of norm on a certain maximal torus inside $\mathrm{G}$,
and in general it is difficult to compute.

In the case of the action of $\GLd$ on $\Rep(Q, \tuple{d})$,
this stratification is characterised by an effective description
of the \emph{Harder--Narasimhan filtration} of each representation.

For a stability parameter $\theta$ and some $\tuple{a} \in \mathbb{Z}_{>0}^{Q_0}$,
we define the associated \emph{slope function} as
\begin{equation}
  \mu\colon\mathbb{Z}_{\geq 0}^{Q_0} \setminus \{\tuple{0}\} \to \mathbb{Q} : \tuple{e} \mapsto \frac{\theta(\tuple{e})}{\tuple{a}\cdot\tuple{e}}.
\end{equation}
The standard choice is $\tuple{a} = \tuple{1}$,
for which~$\mathbf{a}\cdot\mathbf{e}=\sum_{i\in Q_0}e_i$.
Different choices of $\tuple{a}$
correspond to different choices of
norm on the maximal torus of $\GLd$,
see~\cite[Theorem~3.8]{MR3871820}.

\begin{definition}
  Let $Q$ and $\tuple{d}$ be a quiver and a dimension vector,
  and let $\mu$ be a slope function.
  A representation $V \in \Rep(Q,\tuple{d})$ is \emph{$\mu$-stable},
  respectively \emph{$\mu$-semistable},
  if every non-zero subrepresentation $W \subset V$ satisfies
  the inequality $\mu(\dimvec(W)) < \mu(\dimvec(V))$,
  respectively~$\mu(\dimvec(W)) \leq \mu(\dimvec(V))$.
\end{definition}
Note that we do not require that $\mu(\dimvec(V)) = 0$.

This more general notion of stability
allows us to associate to each representation a unique
filtration by subrepresentations,
as it is done in \cite[Theorem~2.5]{MR1906875}
and independently \cite[Proposition~2.5]{MR1974891}.

\begin{proposition}
  Let $Q$ and $\tuple{d}$ be a quiver and a dimension vector,
  and let $\mu$ be a slope function.
  Every representation $V \in \Rep(Q,\tuple{d})$ admits a unique filtration
  \begin{equation}
    0 = V^{0} \subset V^{1}\subset\dots\subset V^{\ell} = V,
  \end{equation}
  which has the property that
  \begin{itemize}
    \item successive quotients $W^{i} \colonequals V^{i}/V^{i-1}$
      have strictly decreasing slopes, and
    \item each $W^{i}$ is $\mu$-semistable.
  \end{itemize}
\end{proposition}
We call this filtration the \emph{Harder--Narasimhan filtration}
of $V$ associated to $\mu$.
The sequence of dimension vectors $(\dimvec(W^{1}), \dots, \dimvec(W^{\ell}))$
is called the \emph{Harder--Narasimhan type} of $V$, and we denote it
by~$\tuple{d}^* = (\tuple{d}^1,\dots,\tuple{d}^{\ell})$.

We see that a representation $V$ is $\mu$-semistable
if and only if its Harder--Narasimhan type is the trivial one, i.e., $(\dimvec(V))$.

One can show,
cf.~\cite[Proposition~3.4]{MR3871820},
that the Hesselink stratification of the unstable locus
is characterised by Harder--Narasimhan types:
two representations in $\Rep(Q,\tuple{d})$ belong to the same stratum
if and only if their Harder--Narasimhan types are equal.
To enumerate all the Hesselink strata it suffices then
to enumerate all the possible Harder--Narasimhan types,
which can be found through combinatorial arguments.

We illustrate these notions with \quivertools.

\iftoggle{sagecode}{\iftrue}{\iffalse}
\begin{sagesnippet}
|\sage| Q = KroneckerQuiver(4); d = (2, 3); theta = (3, -2);
|\sage| M = QuiverModuliSpace(Q, d, theta)
|\sage| M.all_harder_narasimhan_types()
(((2, 3),),
 ((1, 1), (1, 2)),
 ((2, 2), (0, 1)),
 ((2, 1), (0, 2)),
 ((1, 0), (1, 3)),
 ((1, 0), (1, 2), (0, 1)),
 ((1, 0), (1, 1), (0, 2)),
 ((2, 0), (0, 3)))
\end{sagesnippet}
\fi
\iftoggle{juliacode}{\iftrue}{\iffalse}
\begin{minted}{jlcon}
julia> Q, d, theta = kronecker_quiver(4), [2, 3], [3, -2];

julia> M = QuiverModuliSpace(Q, d, theta);

julia> all_hn_types(M)
8-element Vector{HNType}:
 [[2, 3]]
 [[1, 1], [1, 2]]
 [[2, 2], [0, 1]]
 [[2, 1], [0, 2]]
 [[1, 0], [1, 3]]
 [[1, 0], [1, 2], [0, 1]]
 [[1, 0], [1, 1], [0, 2]]
 [[2, 0], [0, 3]]
\end{minted}
\fi

\paragraph{Betti numbers of quiver moduli}
One direct consequence of this enumeration is the ability to compute
Betti numbers for quiver moduli.
This is done in \cite[\S6]{MR1974891},
by translating the method originally introduced by Harder and Narasimhan in \cite{MR0364254}
to compute the Betti numbers of moduli spaces of vector bundles on curves.
We refer to \cite{MR1974891} for more details.
We compute the Betti numbers of our running example using \quivertools.

\iftoggle{sagecode}{\iftrue}{\iffalse}
\begin{sagesnippet}
|\sage| Q = KroneckerQuiver(4); d = (2, 3); theta = (3, -2);
|\sage| M = QuiverModuliSpace(Q, d, theta)
|\sage| M.betti_numbers()
[1, 0, 1, 0, 3, 0, 4, 0, 7, 0, 8, 0, 10, 0, 8, 0, 7, 0, 4, 0, 3, 0, 1, 0, 1]
\end{sagesnippet}
\fi
\iftoggle{juliacode}{\iftrue}{\iffalse}
\begin{minted}{jlcon}
julia> Q, d, theta = kronecker_quiver(4), [2, 3], [3, -2];

julia> M = QuiverModuliSpace(Q, d, theta);

julia> print(betti_numbers(M))
[1, 0, 1, 0, 3, 0, 4, 0, 7, 0, 8, 0, 10, 0, 8, 0, 7, 0, 4, 0, 3, 0, 1, 0, 1]
\end{minted}
\fi

\section{Teleman quantization}
\label{section:teleman-quantization}
In \cref{section:harder-narasimhan} we introduced the Hesselink
stratification of a variety $X$
under the action of a reductive group $\mathrm{G}$.
It is used in \cite{MR3327537}
to obtain \emph{Teleman quantization},
which makes it possible to study the cohomology
of coherent sheaves on $X \gitquot\mathrm{G}$.
For quiver moduli problems,
Teleman quantization can be used effectively
thanks to computations performed in \cite{2311.17003}.
The authors already used this tool in \cite{2311.17003, MR4894756}
to prove several results about rigidity,
existence of partial tilting objects and
semiorthogonal decompositions.

\paragraph{Setup}
We will introduce the necessary notation in the general setting of GIT,
and subsequently specialize it to quiver moduli.

In the Hesselink stratification, each stratum $S$
of the unstable locus $X \setminus X^{\sstable}$ is
identified by a one-parameter subgroup $\lambda_{S}$ of $\mathrm{G}$,
that is in a certain sense ``the most responsible'' subgroup
for the instability of the points in $S$.
Each such $\lambda_{S}$ acts on $\det({\normal^{\vee}_{S/X}})|_{S^{\lambda_{S}}}$,
with an integer weight that we denote by $\eta_{S}$.

Let $F$ be a linearised coherent sheaf on $X$
such that $F|_{X^{\sstable}}$
descends to $\mathcal{F}$ on the quotient~$X^{\sstable}\gitquot\mathrm{G}$.
For each stratum $S$,
denote the set of $\lambda_{S}$-weights of $F|_{S^{\lambda_{S}}}$ by $\operatorname{W}(F, S)$.
Teleman quantization then states that,
if all the $\lambda_{S}$-weights of $F$ on each stratum
are strictly smaller than $\eta_{S}$,
then there exists an isomorphism
\begin{equation}
  \HH^{\bullet}(X^{\sstable}, F|_{X^{\sstable}})^{\mathrm{G}} \cong \HH^{\bullet}(X, F)^{\mathrm{G}}.
\end{equation}
Once again, this result is hard to apply effectively in general, because
describing the Hesselink stratification is a priori difficult.
In the case of quiver moduli however,
the latter is characterised by the Harder--Narasimhan types
(see \cite[Theorem~3.11]{2311.17003}),
so it can be described effectively.

\paragraph{For quiver moduli}
The variety of interest is $\Rep(Q,\tuple{d})$,
with the action of $\GLd$,
and we index all the terms introduced above by
their Harder--Narasimhan types $\tuple{d}^*$.
For each stratum $S_{\tuple{d}^*}$,
the corresponding one-parameter subgroup of $\GLd$
is denoted by~$\lambda_{\tuple{d}^*}$.
For every $\GLd$-equivariant coherent sheaf $F$ on $\Rep(Q,\tuple{d})$,
we denote the set of $\lambda_{\tuple{d}^*}$-weights of $F|_{{S_{\tuple{d^*}}}^{\lambda_{\tuple{d}^*}}}$
by $\operatorname{W}(F, \tuple{d}^*)$.
The weight of the natural action of $\lambda_{\tuple{d}^*}$
on $\det(\normal^{\vee}_{S_{\tuple{d}^*}/\Rep(Q,\tuple{d}^*)})|_{{S_{\tuple{d^*}}}^{\lambda_{\tuple{d}^*}}}$
is denoted by $\eta_{\tuple{d}^*}$.

In the context of quiver moduli problems, Teleman quantization is stated as follows.
\begin{theorem}
  With the previous notation,
  let $F$ be a $\GLd$-equivariant coherent sheaf that descends to $\mathcal{F}$.
  If on every stratum $S_{\tuple{d}^*}$ we have
  \begin{equation}
    \label{equation:teleman-inequality}
    \max \operatorname{W}(F, \tuple{d}^*) < \eta_{\tuple{d}^*},
  \end{equation}
  then for every $k \in \mathbb{Z}$
  we have
  \begin{equation}
    \label{equation:teleman-isomorphisms}
    \HH^{k}(\Rep(Q,\tuple{d}),F)^{\GLd} \cong
    \HH^{k}(\Rep^{\theta\semistable}(Q,\tuple{d}),F|_{\Rep^{\theta\semistable}(Q,\tuple{d})})^{\GLd} \cong
    \HH^{k}(\modulispace^{\theta\semistable}(Q,\tuple{d}), \mathcal{F}).
  \end{equation}
\end{theorem}

In particular, since $\Rep(Q,\tuple{d})$ is affine,
we have vanishing of higher cohomology.
\begin{corollary}
  Under the previous assumptions, we have for all $k \geq 1$ the vanishing
  \begin{equation}
    \HH^{k}(\modulispace^{\theta\semistable}(Q,\tuple{d}), \mathcal{F}) = 0.
  \end{equation}
\end{corollary}

The effective computations of the numbers~$\eta_{\tuple{d}^*}$
are performed in \cite[Corollary~3.18]{2311.17003},
where these are shown to be expressed
in terms of the Euler form of $Q$.
\quivertools implements these computations,
making it possible to verify \eqref{equation:teleman-inequality}.
We illustrate this in \quivertools using our running example.

\iftoggle{sagecode}{\iftrue}{\iffalse}
\begin{sagesnippet}
|\sage| Q = KroneckerQuiver(4); d = (2, 3); theta = (3, -2);
|\sage| M = QuiverModuliSpace(Q, d, theta);
|\sage| hn_types = M.all_harder_narasimhan_types(proper=True)
|\sage| table(
|\sagedots|    [["Harder-Narasimhan type hn", "eta_hn"]] + \
|\sagedots|    [[hn, M.teleman_bound(hn)] for hn in hn_types]
|\sagedots| )
  Harder-Narasimhan|{\tt type}| hn   eta_hn
  ((1, 1), (1, 2))            25
  ((2, 2), (0, 1))            30
  ((2, 1), (0, 2))            70
  ((1, 0), (1, 3))            55
  ((1, 0), (1, 2), (0, 1))    140
  ((1, 0), (1, 1), (0, 2))    125
  ((2, 0), (0, 3))            120
\end{sagesnippet}
\fi
\iftoggle{juliacode}{\iftrue}{\iffalse}
\begin{minted}{jlcon}
julia> Q, d, theta = kronecker_quiver(4), [2, 3], [3, -2];

julia> M = QuiverModuliSpace(Q, d, theta);

julia> hn_types = all_hn_types(M; unstable=true);

julia> weights = hcat(map(hn -> [hn, [teleman_bounds(M)[hn]]], hn_types)...);

julia> permutedims(hcat(["Harder-Narasimhan type hn", "eta_hn"], weights), (2, 1))
8×2 Matrix{Any}:
 "Harder-Narasimhan type hn"  "eta_hn"
 [[1, 1], [1, 2]]             [25]
 [[2, 2], [0, 1]]             [30]
 [[2, 1], [0, 2]]             [70]
 [[1, 0], [1, 3]]             [55]
 [[1, 0], [1, 2], [0, 1]]     [140]
 [[1, 0], [1, 1], [0, 2]]     [125]
 [[2, 0], [0, 3]]             [120]
\end{minted}
\fi

\paragraph{Universal families on quiver moduli}
An application of the Teleman quantization theorem
are vanishing results for the cohomology
of interesting vector bundles on quiver moduli.

Let $Q$ be an acyclic quiver and $\tuple{d}$ be a coprime dimension vector
for which stable representations exist.
For a sufficiently general choice of stability parameter $\theta$,
the moduli space $\modulispace^{\theta\semistable}(Q, \tuple{d})$
is a smooth projective variety admitting
a universal family~$\mathcal{U} = \bigoplus_{i\in Q_0}\mathcal{U}_{i}$.
This universal family has the property that
for a closed point in the moduli space,
the fiber of the bundle~$\mathcal{U}$
is the representation of~$Q$ corresponding to that point
of the moduli space.
Its summands appear in the
four-term exact sequence
\begin{equation}
  \label{equation:4-term-global-sequence}
  0\to
  \mathcal{O}_{\modulispace^{\theta\semistable}(Q,\tuple{d})}\to
  \bigoplus_{i \in Q_0}\mathcal{U}^{\vee}_{i}\otimes\mathcal{U}_{i}\to
  \bigoplus_{\alpha \in Q_1}\mathcal{U}^{\vee}_{\source(\alpha)}\otimes\mathcal{U}_{\target(\alpha)}\to
  \mathrm{T}_{\modulispace^{\theta\semistable}(Q,\tuple{d})}\to
  0.
\end{equation}
This is the globalisation of the exact sequence in \eqref{equation:4-term-sequence},
see \cite[Lemma~4.2]{MR4770368} for its construction.

The Teleman weights of~$\mathcal{U}_{i}$
and~$\mathcal{U}_i^\vee\otimes\mathcal{U}_j$
are computed in \cite[Lemma~3.19 and Proposition~3.20]{2311.17003}.

\quivertools implements the computation of these weights,
as well as the computation of weights
of their endomorphism bundles \cite[Proposition~3.20]{2311.17003},
and, under mild conditions,
of the canonical line bundle $\omega_{\modulispace^{\theta\semistable}(Q,\tuple{d})}$,
which is described in \cite[Proposition~4.7]{MR4894756}.
We illustrate this using in \quivertools using our running example.

\iftoggle{sagecode}{\iftrue}{\iffalse}
\begin{sagesnippet}
|\sage| Q = KroneckerQuiver(4); d = (2, 3); theta = (3, -2);
|\sage| M = QuiverModuliSpace(Q, d, theta);
|\sage| U0 = M.weights_universal_bundle(0)
|\sage| U1 = M.weights_universal_bundle(1)
|\sage| omega = M.weights_canonical()
|\sage| hn_types = M.all_harder_narasimhan_types(proper=True)
|\sage| table(
|\sagedots|     [["Harder-Narasimhan type hn", "W(U0, hn)", "W(U1, hn)", "W(omega, hn)"]] + \
|\sagedots|     [[hn, U0[hn], U1[hn], omega[hn]] for hn in hn_types]
|\sagedots| )
  Harder-Narasimhan|{\tt type}| hn   W(U0, hn)   W(U1, hn)      W(omega, hn)
  ((1, 1), (1, 2))            [5, 0]      [5, 0, 0]      [20]
  ((2, 2), (0, 1))            [5, 5]      [5, 5, 0]      [40]
  ((2, 1), (0, 2))            [10, 10]    [10, 5, 5]     [80]
  ((1, 0), (1, 3))            [10, 5]     [5, 5, 5]      [60]
  ((1, 0), (1, 2), (0, 1))    [25, 15]    [15, 15, 10]   [160]
  ((1, 0), (1, 1), (0, 2))    [20, 15]    [15, 10, 10]   [140]
  ((2, 0), (0, 3))            [15, 15]    [10, 10, 10]   [120]
\end{sagesnippet}
\fi
\iftoggle{juliacode}{\iftrue}{\iffalse}
\begin{minted}{jlcon}
julia> Q, d, theta = kronecker_quiver(4), [2, 3], [3, -2];

julia> M = QuiverModuliSpace(Q, d, theta);

julia> hn_types = all_hn_types(M; unstable=true);

julia> U0, U1, omega = universal_bundle(M, 1), universal_bundle(M, 2), canonical_bundle(M);

julia> weights = hcat(
            map(
                hn -> [hn, map(bundle -> teleman_weights(bundle)[hn], [U0, U1, omega])...],
                hn_types
            )...
        );

julia> permutedims(
            hcat(
                ["Harder-Narasimhan type hn", "W(U0, hn)", "W(U1, hn)", "W(omega, hn)"],
                weights
            ),
            (2,1)
        )
8×4 Matrix{Any}:
 "Harder-Narasimhan type hn"  "W(U0, hn)"  "W(U1, hn)"   "W(omega, hn)"
 [[1, 1], [1, 2]]             [5, 0]       [5, 0, 0]     [20]
 [[2, 2], [0, 1]]             [5, 5]       [5, 5, 0]     [40]
 [[2, 1], [0, 2]]             [10, 10]     [10, 5, 5]    [80]
 [[1, 0], [1, 3]]             [10, 5]      [5, 5, 5]     [60]
 [[1, 0], [1, 2], [0, 1]]     [25, 15]     [15, 15, 10]  [160]
 [[1, 0], [1, 1], [0, 2]]     [20, 15]     [15, 10, 10]  [140]
 [[2, 0], [0, 3]]             [15, 15]     [10, 10, 10]  [120]
\end{minted}
\fi

\section{Intersection theory}
\label{section:intersection-theory}
When quiver moduli are smooth projective varieties
admitting a universal family,
an effective presentation of their Chow rings,
together with an expression for the Todd class and point class,
is given in \cite{MR3318266,MR4770368}.

The presentation of the Chow ring
has the Chern classes of the summands of
the universal family~$\bigoplus_{i\in Q_0}\mathcal{U}_i$ as generators,
and depends on the choice of linearisation $\chi$.
This description makes it possible to compute Euler characteristics of sheaves
admitting an expression in the Chern classes of the~$\mathcal{U}_i$,
or compute degrees of line bundles,
as in \cite{MR4770368,2412.15390v1}.

Chow ring calculations are numerical
and only determine the Euler characteristic,
but these results can be used together with
vanishing results from \cref{section:teleman-quantization}
to know the number of global sections of certain vector bundles.
See also \cite{MR4883104}
for canonical identifications of global sections of~$\mathcal{U}_i^\vee\otimes\mathcal{U}_j$.

We illustrate these features in \quivertools using our running example.

\iftoggle{sagecode}{\iftrue}{\iffalse}
\begin{sagesnippet}
|\sage| Q = KroneckerQuiver(4); d = (2, 3); theta = (3, -2);
|\sage| M = QuiverModuliSpace(Q, d, theta);
|\sage| chi = (-1, 1);
|\sage| CH = M.chow_ring(chi=chi)
|\sage| CH.gens()
(x1_1bar, x0_2bar, x1_1bar, x1_2bar, x1_3bar)
|\sage| td = M.todd_class(chi=chi);
|\sage| pt = M.point_class(chi=chi);
|\sage| # Computing the Hilbert series of M
|\sage| H = M.chern_character_line_bundle(theta)
|\sage| [M.integral(H ** i) for i in range(M.dimension())]
[1,
 126,
 4032,
 59268,
 531839,
 3395882,
 16907632,
 69626910,
 246885947,
 775675824,
 2205490144,
 5766791394]
|\sage| # Computing the degree of -K_M^dim(M)
|\sage| anticanonical = H ** M.index();
|\sage| (anticanonical ** M.dimension()).lift().homogeneous_components()[M.dimension()] / pt
1996824248320
\end{sagesnippet}
\fi
\iftoggle{juliacode}{\iftrue}{\iffalse}
\begin{minted}{jlcon}
julia> Q = kronecker_quiver(4); d = [2, 3]; theta = canonical_stability(Q, d);

julia> M = QuiverModuliSpace(Q, d, theta);

julia> CH = chow_ring(M; chi=[-1, 1])
Singular polynomial quotient ring (QQ),(x11,x12,x21,x22,x23),(dp(5),C)

julia> td = todd_class(M); pt = point_class(M);

julia> H = line_bundle(M, Int.(1/index(M) * theta))
Bundle of rank 1

julia> [integral(H^i) for i in 0:dimension(M)-1] # Computing the Hilbert series of M
12-element Vector{Singular.n_Q}:
 1
 126
 4032
 59268
 531839
 3395882
 16907632
 69626910
 246885947
 775675824
 2205490144
 5766791394

julia> degree(dual(canonical_bundle(M))) # Computing the degree of -K_M^dim(M)
 1996824248320
\end{minted}
\fi

\printbibliography
\vfill
\emph{Pieter Belmans}, \url{p.belmans@uu.lu} \\
Mathematical Institute, Utrecht University, Budapestlaan 6, 3584CD Utrecht, Netherlands

\emph{Hans Franzen}, \url{hans.franzen.math@gmail.com}

\emph{Gianni Petrella}, \url{gianni.petrella@uni.lu} \\
Department of Mathematics, Universit\'e de Luxembourg, 2, Avenue de l'Université, 4365 Esch-sur-Alzette, Luxembourg

\end{document}

%% file: preamble.tex
\usepackage{fontenc}
\usepackage{etoolbox}
\usepackage{graphicx, float}
\usepackage[a4paper, total={6in, 9.5in}]{geometry}
\usepackage{parskip}

\usepackage{colonequals}
\usepackage{booktabs}
\usepackage[fleqn,leqno]{amsmath}
\usepackage{amsthm, amsfonts, amssymb}

\usepackage{hyperref, bookmark}
\hypersetup{hypertexnames = false, bookmarksdepth = 2, bookmarksopen = true, colorlinks, linkcolor = black, citecolor = black, urlcolor = black, pdfstartview={XYZ null null 1}}
\usepackage[capitalise]{cleveref}

\usepackage{xcolor}
\definecolor{LightGray}{gray}{0.95}

\usepackage[frozencache=true,cachedir=minted-cache]{minted}
\setminted{
frame=lines,
bgcolor=LightGray,
baselinestretch=1.0,
fontsize=\footnotesize,
framesep=1mm,
}

\usepackage{tikz, tikz-cd}
\usetikzlibrary{calc}
\usetikzlibrary{arrows}
\usetikzlibrary{matrix}
\usetikzlibrary{positioning}

\usepackage[T1]{fontenc}
\usepackage{libertine}
\usepackage[libertine]{newtxmath}
\usepackage[scaled=0.8]{FiraMono}
\usepackage[backend=biber, datamodel=mrnumber, maxbibnames=99, sortcites]{biblatex}
\usepackage{gitinfo2}

\usepackage{thmtools}
\declaretheoremstyle[
  spaceabove = 3pt,
  spacebelow = 3pt,
  bodyfont = \itshape,
]{first}
\declaretheoremstyle[
  spaceabove = 3pt,
  spacebelow = 3pt,
]{second}
\declaretheorem[numberwithin=section, style=first]{theorem}

\declaretheorem[sibling=theorem, style=first]{corollary}

\declaretheorem[sibling=theorem, style=first]{proposition}

\declaretheorem[sibling=theorem, style=second]{definition}

\Crefname{assumption}{Assumption}{Assumptions}
\Crefname{convention}{Convention}{Conventions}
\Crefname{setup}{Setup}{Setups}

\crefname{alphatheorem}{Theorem}{Theorems}
\crefname{alphaconjecture}{Conjecture}{Conjectures}
\crefname{alphaquestion}{Question}{Questions}

\usepackage{xspace}

\newcommand\dash{\nobreakdash-\hspace{0pt}}

\DeclareMathOperator{\Hom}{Hom}
\DeclareMathOperator{\Ext}{Ext}
\DeclareMathOperator{\ext}{ext}
\DeclareMathOperator{\GL}{GL}
\DeclareMathOperator{\GLd}{GL(\tuple{d})}
\DeclareMathOperator{\Rep}{R}
\DeclareMathOperator{\modulispace}{M}

\DeclareMathOperator{\HH}{H}

\DeclareMathOperator{\normal}{N}

\DeclareMathOperator{\sheafHom}{\mathscr{H}\kern -2.5pt\mathit{om}}

\newcommand{\stable}{\hyphen\mathrm{st}}
\newcommand{\semistable}{\hyphen\mathrm{sst}}
\newcommand{\sstable}{\mathrm{sst}}
\newcommand{\tuple}[1]{\ensuremath{\mathbf{#1}}}

\newcommand\field{\ensuremath{\mathbf{k}}}

\DeclareMathOperator{\target}{t}
\DeclareMathOperator{\source}{s}
\DeclareMathOperator{\Spec}{Spec}
\DeclareMathOperator{\Mat}{Mat}

\newcommand{\quivertools}{{\tt{QuiverTools}}\xspace}
\newcommand{\dimvec}{\underline{\dim}}
\newcommand{\gitquot}{/\kern-0.2em /}

\newtoggle{sagecode}
\newtoggle{juliacode}
\toggletrue{sagecode}
\togglefalse{juliacode}

\newenvironment{sagesnippet}{\VerbatimEnvironment
  \begin{minted}[escapeinside=||]{sage}}
 {\end{minted}}

\newcommand{\sage}{{\color{blue}\texttt{sage:}}}
\newcommand{\sagedots}{{\color{blue}\texttt{....:}}}